\documentclass[a4paper,12pt]{amsart}

\usepackage[all]{xy}
\usepackage[dvips]{graphicx}
\usepackage{psfrag}

\newtheorem{teo}{Theorem}[section]
\newtheorem{defi}[teo]{Definition}
\newtheorem{lema}[teo]{Lemma}
\newtheorem{obs}[teo]{Remark}
\newtheorem{cor}[teo]{Corollary}

\newcommand{\bdefi}{\begin{defi}}
\newcommand{\edefi}{\end{defi}}
\newcommand{\bteo}{\begin{teo}}
\newcommand{\eteo}{\end{teo}}
\newcommand{\blem}{\begin{lema}}
\newcommand{\elem}{\end{lema}}
\newcommand{\bobs}{\begin{obs}}
\newcommand{\eobs}{\end{obs}}
\newcommand{\bcor}{\begin{cor}}
\newcommand{\ecor}{\end{cor}}

\newcommand{\D}{\mathfrak{D}}
\newcommand{\fimdem}{\phantom. \hfill \square}
\newcommand{\Z}{\mathbb{Z}}
\newcommand{\dem}{{\bf Proof: }}
\newcommand{\ov}{\overline}
\newcommand{\tl}{\widetilde}
\newcommand{\sbs}{\subset}

\begin{document}  

\title{Coincidence classes in nonorientable manifolds}
\author{Daniel Vendr\'uscolo} 
\email{vendrus@picard.ups-tlse.fr}

\begin{abstract}
In this article we studied Nielsen coincidence theory for maps between
manifolds of same dimension without hypotheses on orientation. We use the
definition of semi-index of a class, we review the definition of defective
classes and study the appearance of defective root classes. We proof a
semi-index product formula type for lifting maps and we presented conditions
such that defective coincidence classes are the only essencial classes.
\end{abstract}


\maketitle

\section{Introduction}     

Nielsen coincidence theory was extended (\cite{doje} and \cite{je2}) to
maps between nonorientable topological manifolds using the notion of
$semi-index$ (a non negative integer) for a coincidence set. 

We consider maps $f,g:M\to N$ between manifolds without boundary of the same
dimension $n$, we define $h=(f,g):M\to N\times N$, then using microbundles (see
\cite{je2} for details) we can suppose that $h$ is in a transverse position. 

Let $w$ be a path satisfying the Nielsen relation between $x,y\in
Coin(f,g)$. We choose a local orientation $\gamma_0$ of $M$ in $x$ and denote
by $\gamma_t$ the translation of $\gamma_0$ along $w(t)$.

\bdefi \cite[1.2]{je2} We will say that two points $x,y\in
Coin(f,g)$ are $R$-related $(xRy)$ if and only if there is a path $w$
establishing the Nielsen relation between them such that the translation of
the orientation $h_*\gamma_0$ along a path in the diagonal $\Delta(N)\sbs
N\times N$ homotopic to $h\circ w$ in $N\times N$ is opposite to
$h_*\gamma_1$. In this case the path $w$ is called graph-orientation-reversing.
\edefi

Since $(f,g)$ is transverse, $Coin(f,g)$ is finite. Let $A\sbs Coin(f,g)$,
then $A$ can be represented as
$A=\{a_1,a_2,\cdots,a_s;b_1,c_1,\cdots,b_k,c_k\}$ where $b_iRc_i$ for any $i$
and $a_iRa_j$ for no $i\not= j$. The elements $\{a_i\}_i$ of this
decomposition are called {\it free}.

\bdefi In the above situation the $semi-index$ of the pair $(f,g)$ in
$A=\{a_1,a_2,\cdots,a_s;b_1,c_1,\cdots,b_k,c_k\}$ 
(denoted $|ind|(f,g;A)$) is the number of free elements of this decomposition
of $A$.  
$$|ind|(f,g;A)=s.$$
\edefi

In \cite{doje} and \cite{je2} we can find the proof of the fact that this
definition does not depend on the decomposition of $A$. Moreover, if $U\sbs M$
is an open
subset we can extend this definition to a semi-index of a pair on the
subset $U$ ($|ind|(f,g;U)$).

\bdefi A coincidence class $C$ of a transverse pair $(f,g)$ is essential if  
$|ind|(f,g;A)\not=0.$
\edefi

In \cite{je} the authors studied when a coincidence point $x\in Coin(f,g)$ satisfy
$xRx$. These points can only appear in the nonorientable case, they will be
called {\it self-reducing} points. They
represent a new situation (see \cite[Example 2.4]{je}) that never occurs in
the orientable case or in the fixed point context.

\bdefi (\cite[2.1]{je}) Let $x\in Coin(f,g)$ and let $H\sbs \pi_1(M)$, $H'\sbs
\pi_1(N)$ denote the subgroups of orientation-preserving elements. We define:
$$Coin(f_{\#},g_{\#})_x=\{\alpha\in \pi_1(M,x)\ |\
f_{\#}(\alpha)=g_{\#}(\alpha)\},$$
$$Coin^+(f_{\#},g_{\#})_x=Coin(f_{\#},g_{\#})_x\cap H.$$
\edefi

\blem \label{defloop} (\cite[2.2]{je}) Let $f,g:M\to N$ be transverse and
$x\in Coin(f,g)$. Then $xRx$ if and only if $Coin^+(f_{\#},g_{\#})_x\not=
Coin(f_{\#},g_{\#})_x\cap f_{\#}^{-1}(H')$ (in other words, if there exists a
loop $\alpha$ based at $x$ such that $f\circ\alpha \sim g\circ\alpha$ and
exactly one of the loops $\alpha$ or $f\circ\alpha$ is orientation-preserving).
\elem

\bdefi A coincidence class $A$ is called defective if $A$ contains a
self-reducing point.
\edefi

\blem (\cite[2.3]{je})  \label{inddef} If a Nielsen class $A$ contains a
self-reducing point (i.e. $A$ is defective) then any two points in this class
are $R$-related, and thus:
$$|ind|(f,g;A)=\left\{\begin{array} {cc}
0 & \mbox{if $\#A$ is even}; \\ 
1 &  \mbox{if $\#A$ is odd}.
\end{array} \right.$$
\elem

\section{The root case}

In \cite{bs} we can find a different approach to extending the Nielsen root
theory to the nonorientable case. Using the concept of {\it
orientation-true}\footnote{$f$ is orientation-true if for each loop 
$\alpha\in \pi_1(M)$, $f(\alpha)$ is orientation-preserving if and only if $\alpha$ is
orientation-preserving.} map they classified maps between manifolds of the
same dimension in three
types (see also \cite{olum} and \cite{sk}):

\bdefi (\cite[2.1]{bs}) Let $f:M\to N$ be a map of manifolds. Then three types
of maps are defined as follows.
\begin{itemize}
\item[(1)] Type I: $f$ is orientation-true.
\item[(2)] Type II: $f$ is not orientation-true but does not map an
  orientation-reversing loop in $M$ to a contractible loop in $N$.
\item[(3)] Type III: $f$ maps an orientation-reversing loop in $M$ to a
  contractible loop in $N$. 
\end{itemize}
Further, a map $f$ is defined to be orientable if it is of Type I or II, and
nonorientable otherwise.
\edefi

For orientable maps they describe an {\it Orientation Procedure}
(\cite[2.6]{bs}) for
root classes, using local degree with coefficients in $\Z$. For maps of Type III
the same procedure is only possible with coefficients in $\Z_2$.
They then defined the {\it multiplicity} of a root class, that is an integer
for orientable maps and an element of $\Z_2$ for maps of Type III

Now if we consider the root classes for a map $f$ to be the coincidence classes
of the pair $(f,c)$ where $c$ is the constant map, we have:

\bteo \label{rootdef} Let $f:M\to N$ be a map of manifolds of the same
dimension. 
\begin{itemize}
\item[(i)] If $f$ is orientable, no root class of $f$ is defective.
\item[(ii)] If $f$ is of Type III, then all root classes of $f$ are defective.
\end{itemize}
\eteo

\dem By lemma~\ref{defloop}, a coincidence class $C$ of the pair $(f,c)$ is
defective if and only if there exists a point $x\in C$ and a loop $\alpha$ in $x$ such
that $f\circ\alpha \sim 1$ and $\alpha$ is orientation-reversing. This proves
the first statement. 

Now let $f$ be a Type III map. Then there exists a loop $\alpha\in
\pi_1(M,x_0)$ such that $\alpha$ is orientation-reversing and
$f\circ\alpha\sim 1$. If $x$ is a root of $f$ choosing a path $\beta$ from
$x$ to $x_0$ we have that $\gamma=\beta^{-1}\alpha\beta$ is a loop in $x$ such
that $\gamma$ is orientation-reversing and $f\circ\gamma\sim 1$. Then all
roots of $f$ are self-reducing points. $\fimdem$

In fact \cite[Lemma~4.1]{bs} shows the equality between the multiplicity of a
root class and its semi-index.

\bteo \label{raizes} Let $M$ and $N$ be manifolds of the same dimension such
that $M$ is nonorientable and $N$ is orientable. If $f:M \to N$ is a map then
all essential root classes of $f$ are defective.  \eteo

\dem There are no orientation-true maps from a nonorientable to an orientable
manifold. If $f$ is a Type II map then by \cite[Lemma~3.10]{bs} $deg(f)=0$ and
$f$ has no essential root classes. The result follows by
proposition~\ref{rootdef}. $\fimdem$ 

Using the ideas of Proposition~\ref{rootdef} we can also state:

\blem \label{loopcen} Let $f,g:M\to N$ be two maps between manifolds of the
same dimension. If there exist a coincidence point $x_0$ and a
graph-orientation-reverse loop $\alpha$ in $x_0$ such that $f(\alpha)$ is in
the center of $\pi_1(N)$, then all coincidence points of the pair $(f,g)$ are
self-reducing points.
\elem

\dem We can suppose that all coincidences of the pair $(f,g)$ have image on the
same point in $N$. If $x_1\in Coin(f,g)$ we choose a path $\beta$ from $x_1$ to $x_0$
and take the loop $\gamma=\beta^{-1}\circ\alpha\circ\beta$ at $x_1$. Since
$\alpha$ belongs to the center of $\pi_1(N)$,  $\gamma$
is a graph-orientation-reverse loop at $x_1$. $\fimdem$

\bcor \label{imacen} Let $f,g:M\to N$ be two maps between manifolds of the
same dimension such that  $f_{\#}(\pi_1(M))$ is contained in the center of
$\pi_1(N)$. If $(f,g)$ has a defective class then all classes of $(f,g)$ are
defective. 
\ecor
 
\section{Covering maps}

Let $M$ and $N$ be compact, closed manifolds of the same dimension; 
$f,g:M\to N$ be two maps such that $Coin(f,g)$ is finite. and
 $p:\tl M\to M$ and $q:\tl N \to N$ 
be finite coverings such that there exist lifts $\tl f,\tl g:\tl M\to \tl N$
of the pair $f,g$.
$$\xymatrix{
\tl{M} \ar[r]^-{\tl{f}} _-{\tl g} \ar[d]_-{p} & \tl N \ar[d]^-{q} \\
M \ar[r]^-{f} \ar[r]_-{g} & {N}}
$$

In this situation the groups of Deck transformations of the lifts $\tl
M$ and $\tl N$ can be described by:
$$\D(\tl M)=\frac{\pi_1(M)}{p_{\#}(\pi_1(\tl M))};  \ \ \ 
\D(\tl N)=\frac{\pi_1(N)}{q_{\#}(\pi_1(\tl N))}$$

Choosing a point ${x_0}\in Coin(f,g)$ we define $\tl f_{*,x_0},
\tl g_{*,x_0} :\D(\tl M)\to \D(\tl N)$ such that, for each $\alpha\in \D(\tl
M)$ we take 
$\ov \alpha\in \pi_1(M, {x_0})$ such that $\D(\tl M)
\ni [\ov \alpha]=\alpha$ and:
$$\tl f_{*,x_0}(\alpha)=[f_{\#}(\ov\alpha)]\in \D(\tl N)$$  
$$\tl g_{*,x_0}(\alpha)=[g_{\#}(\ov\alpha)]\in \D(\tl N)$$

We note that $\tl f_{*,x_0}, \tl g_{*,x_0}$ may depend on ${x_0}$.
Choosing a ${x_0}\in Coin(f,g)$ we have.   
$\forall \tl x \in\tl M$ and $\forall \alpha \in\D(\tl M)$:
$$\tl f(\alpha(\tl x))=\tl f_{{*,x_0}}(\alpha)\circ \tl f(\tl x)$$  
$$\tl g(\alpha(\tl x))=\tl g_{{*,x_0}}(\alpha)\circ \tl g(\tl x)$$  

\blem \label{lifcoin} Let $\tl{x_0}\in Coin(\tl f,\tl g)$ and $\alpha \in
\D(\tl M)$. Then $\alpha(\tl{x_0})\in Coin(\tl f,\tl g)$ if and only if 
$\tl f_{{*,x_0}}(\alpha)=\tl g_{{*,x_0}}(\alpha)$ where $x_0=p(\tl{x_0})$. 
\elem

\dem By the above observation we have $\tl f(\alpha(\tl{x_0}))= \tl
f_{{*,x_0}}(\alpha)\circ \tl f(\tl{x_0})$ and  $\tl g(\alpha(\tl{x_0}))= \tl
g_{{*,x_0}}(\alpha)\circ \tl g(\tl{x_0})$. This proves the lemma.$\fimdem$

\bcor \label{numcop} Let $\tl{x_0}\in Coin(\tl f,\tl g)$ and
$x_0=p(\tl{x_0})$. Then $p^{-1}(x_0)\cap Coin(\tl f,\tl g)$ have exactly 
$\#Coin(\tl f_{{*,x_0}},\tl g_{{*,x_0}})$ elements.$\fimdem$
\ecor

\blem \label{samcla} Let $\tl{x_0}$ and $\tl{x}'_0$ be two coincidences of the
pair $(\tl f,\tl g)$ such that $p(\tl{x_0})=p(\tl{x}'_0)=x_0$, and let
$\gamma$  the unique element of $\D(\tl M)$
such that $\gamma(\tl{x_0})=\tl{x}'_0$. The points $\tl{x_0}$ and $\tl{x}'_0$
are in the same coincidence class of $(\tl f,\tl g)$ if and only if there
exists $\ov\gamma\in\pi_1(M,x_0)$ such that:
\begin{itemize}
\item $\D(\tl M)\ni [\ov\gamma]=\gamma$;
\item $f_{\#}(\ov\gamma)=g_{\#}(\ov\gamma)$.
\end{itemize}
\elem

\dem ($\Rightarrow$) If $\tl{x_0}$ and $\tl{x}'_0$
are in the same coincidence class of $(\tl f,\tl g)$, there exists a path
$\beta$ from $\tl{x_0}$ to $\tl{x}'_0$ that realizes the Nielsen relation,
({\it i.e.} $\tl f\circ\beta\sim \tl g\circ\beta$).

Take $\ov\gamma=p(\beta) \in \pi_1(M,x_0)$. We can see that
$[\ov\gamma]=\gamma$ and $f\circ\ov\gamma=q\circ\tl f\circ\beta\sim q\circ 
\tl g\circ\beta=g\circ\ov\gamma$ showing that
$f_{\#}(\ov\gamma)=g_{\#}(\ov\gamma)$.

($\Leftarrow$) Take the lift $\tl\gamma$ of $\ov\gamma$ starting at
$\tl{x_0}$. It is a path from $\tl{x_0}$ to $\tl{x}'_0$ that realizes
the Nielsen relation, 
({\it i.e.} $\tl f\circ\tl\gamma\sim \tl g\circ\tl\gamma$).$\fimdem$ 

\bcor \label{lifdef} In lemma~\ref{samcla}, if the points $\tl{x_0}$ and
$\tl{x}'_0$ are in the same coincidence class of $(\tl f,\tl g)$ then
$\tl{x_0}R\tl{x}'_0$ if and only if $sign(\tl f_{*,x_0}(\gamma))\cdot
sign(\gamma)=-1$.  
In this case, $x_0$ is a self-reducing coincidence point.
\ecor

\dem First we note that since $f_{\#}(\ov\gamma)=g_{\#}(\ov\gamma)$ then 
$\tl f_{*,x_0}(\gamma)=\tl g_{*,x_0}(\gamma)$ and we can see that 
$sign(\tl f_{*,x_0}(\gamma))\cdot sign(\gamma)=-1$ if and only if the paths
$\ov\gamma$ and $\tl\gamma$ in the proof of lemma~\ref{samcla} are both graph
orientation-reversing. $\fimdem$

Denoting the natural projection by $j_{x_0}:\pi_1(M, x_0)\to \D(\tl M)$ and
the set $\{\ov\alpha\in\pi_1(M,x_0)\ |\ f_{\#}(\ov\alpha)=g_{\#}(\ov\alpha)\}$
by $Coin(f_{\#},g_{\#})_{x_0}$, we have:

\bcor \label{numsam} If $x_0$ is a coincidence of the pair $(f,g)$ then the
set $p^{-1}(x_0)\cap Coin(\tl f, \tl g)$ can be partitioned in
$\frac{\#Coin(\tl f_{{*,x_0}},\tl g_{{*,x_0}})}
{\#j_{x_0}(Coin(f_{\#},g_{\#})_{x_0})}$ disjoint subsets, each of then with 
$\#j_{x_0}(Coin(f_{\#},g_{\#})_{x_0})$ elements that are all Nielsen related
(therefore contained in the same
coincidence class of the pair $(\tl f, \tl g)$). Moreover, no two points of
different subsets are Nielsen related.$\fimdem$
\ecor

\blem \label{lifpat} Let $x_0, x_1$ be coincidence points in the same
coincidence class of the pair $(f,g)$, $\ov\alpha$ be a path from  $x_0$ to
$x_1$ that realizes the Nielsen relation, $\tl{x_0},\tl{x}'_0$ be coincidence
points of the pair $(\tl f,\tl g)$ such that $p(\tl{x_0})=p(\tl{x}'_0)=x_0$,
and $\gamma$ the unique element of $\D(\tl M)$
such that $\gamma(\tl{x_0})=\tl{x}'_0$. If $\tl\alpha$ and $\tl\alpha'$ are
the two liftings of $\ov\alpha$ starting at $\tl{x_0}$ and $\tl{x}'_0$
respectively then:
\begin{itemize}
\item[(i)] $\tl\alpha(1)$ and $\tl\alpha'(1)$ are coincidence points of the
  pair $(\tl f,\tl g)$;
\item[(ii)] $\tl\alpha(1)$ ($\tl\alpha'(1)$) is in the same coincidence class
  as $\tl{x_0}$ ($\tl{x}'_0$).  
\item[(iii)] $p(\tl\alpha(1))=p(\tl\alpha'(1))=x_1$;
\item[(iv)] $\gamma(\tl\alpha(1))=\tl\alpha'(1)$.
\item[(v)] If $\ov\alpha$ is a graph orientation-reversing-path (in this case
  $x_0Rx_1$) then $\tl\alpha(1)$ and $\tl\alpha'(1)$ are graph
  orientation-reverse-paths (in this case $\tl x_0R\tl x_1$ and $\tl x'_0R\tl
  x'_1$).
\end{itemize}
\elem

\dem (i), (ii) and (iii) are known (we can prove them using covering space
theory). To prove (iv), we take $\ov\gamma\in \pi_1(M, x_0)$ such that $\D(\tl
M)\ni [\ov\gamma]=\gamma$ and we can see that projection of the path
$\ov\alpha\circ \ov\gamma\circ \ov\alpha^{-1}\in \pi_1(M, x_1)$ in $\D(\tl M)$
coincides with $\gamma$.

To prove (v), we use \cite[lemma 2.1, page 77]{doje}.$\fimdem$

\bteo \label{ind} Let $M$ and $N$ be compact, closed manifolds of the same
dimension, $f,g:M\to N$ be two maps such that $Coin(f,g)$ is finite, and
 $p:\tl M\to M$ and $q:\tl N \to N$ be finite coverings such that there
 exist lifts $\tl f,\tl g:\tl M\to \tl N$ of the pair $(f,g)$. If $\tl C$ is a
 coincidence class of the pair $(\tl f,\tl g)$ then $C=p(\tl C)$ is a
 coincidence class of the pair $(f,g)$ and
$$|ind|(\tl f, \tl g; \tl C)= \left\{\begin{array} {cc}
s\cdot k\  (mod\  2) & \mbox{if $C$ is defective}; \\ 
s \cdot k &  \mbox{otherwise},
\end{array} \right.$$
where $s=|ind|(f,g,C)$, $k=\#j(Coin(f_{\#},g_{\#})_{x_0})$ and $x_0\in C$.
\eteo

\dem The fact that $C=p(\tl C)$ is a coincidence class of the pair $(f,g)$ is
known. We choose a point $x_0\in C$. Since $Coin(f,g)$ is finite we can
suppose that a
decomposition 
$C=\{x_1,x_2,\cdots,x_s;c_1,c'_1,c_2,c'_2,\cdots,c_n,c'_n\}$ is such that each
$x_i$ is free, and for all pairs $c_j,c'_j$
we have $c_jRc'_j$.

\begin{figure}[!h]
\centering
\includegraphics[width=9cm]{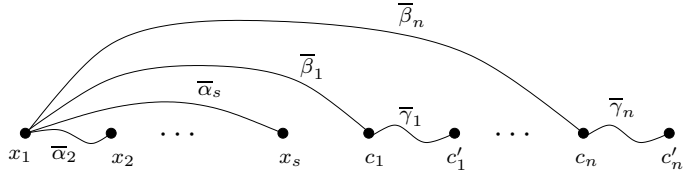}
\caption{The class $C$ and the chosen paths.}
\label{classC}
\end{figure}

Now we choose paths $\{\ov\alpha_i\}_i$, $2\le i\le s$; $\{\ov\beta_j\}_j$ and
$\{\ov\gamma_j\}_j$, $1\le j\le n$ (see figure~\ref{classC}) such that:
\begin{itemize}
\item $\ov\alpha_i$ is a path in $M$ from $x_1$ to $x_i$ that realizes the
  Nielsen relation.
\item $\ov\beta_j$ is a path in $M$ from $x_1$ to $c_j$ that realizes the
  Nielsen relation.
\item $\ov\gamma_j$ is a graph-orientation-reversing path in $M$ from $c_j$ to
  $c'_j$. 
\end{itemize}

\psfrag{a1}{\tiny $\ov\alpha_2$}
\psfrag{b1}{\tiny $\ov\beta_1$}
\psfrag{as}{\tiny $\ov\alpha_s$}
\psfrag{bn}{\tiny $\ov\beta_n$}
\psfrag{g1}{\tiny $\ov\gamma_1$}
\psfrag{gn}{\tiny $\ov\gamma_n$}
\psfrag{cdo}{$\cdots$}
\psfrag{x1}{\tiny $x_1$}
\psfrag{x2}{\tiny $x_2$}
\psfrag{xs}{\tiny $x_s$}
\psfrag{c1}{\tiny $c_1$}
\psfrag{c1l}{\tiny $c'_1$}
\psfrag{cn}{\tiny $c_n$}
\psfrag{cnl}{\tiny $c'_n$}

For each element $\{x^k_1\}$ of $p^{-1}(x_1)\cap Coin(\tl f, \tl g)$ (by
Corollary~\ref{numcop} there were $\#Coin(\tl f_{{*,x_1}},\tl g_{{*,x_1}})$
such elements) we take lifts $\{\tl\alpha^k_i\}_{i,k}$,
$\{\tl\beta^k_j\}_{j,k}$ and $\{\tl\phi^k_j\}_{j,k}$ of the paths
$\{\ov\alpha_i\}_i$, $\{\ov\beta_j\}_j$ 
and $\{\ov\gamma_j\circ\ov\beta_j\}_j$ respectively, starting at $x^k_1$. 

Using Corollary~\ref{numsam} and Lemma~\ref{lifpat} at the point $x_1\in C$ we
 obtain that the set 
 $p^{-1}(C)\cap Coin(\tl f, \tl g)$ is the union of 
$\#Coin(\tl f_{{*,x_1}},\tl g_{{*,x_1}})$ copies of the class $C$ and these
 copies can be divided in
$\frac{\#Coin(\tl f_{{*,x_1}},\tl g_{{*,x_1}})}
{\#j_{x_1}(Coin(f_{\#},g_{\#})_{x_1})}$ disjoint coincidence classes of the
 pair $(\tl f,\tl g)$. Each one of these classes contains 
$\#j_{x_1}(Coin(f_{\#},g_{\#})_{x_1})$ copies of the class $C$.

The result follows by Corollary~\ref{lifdef}.$\fimdem$

\section{Two folded orientable covering}     

Let $M$ and $N$ be manifolds of same dimension with $M$ nonorientable and
$N$ orientable;
$f,g:M\to N$ be two maps such that $Coin(f,g)$ is finite, and $p:\tl M\to M$
be the two-fold orientable covering of $M$. We define $\tl f,\tl g:\tl M \to
N$ by  $\tl f =f\circ p$ and $\tl g = g\circ p$.

$$\xymatrix{
\tl{M} \ar[rd]^-{\tl{f}} _(.4){\tl g} \ar[d]_-{p} & \\
M \ar[r]^-{f} \ar[r]_-{g} & {N}  }
$$

\blem \label{or} Under the above conditions, if $C$ is a coincidence
class of the pair $(f,g)$ then $p^{-1}(C)\subset Coin(\tl f,\tl g)$ is such
that:

\begin{enumerate}

\item $p^{-1}(C)$ can be divided in two disjoint sets $\tl C$ and $\tl
  C'$, of the same cardinality, such that $p(\tl C)=p(\tl C')=C$.

\item If $\tl{x_1}, \tl{x_2}\in \tl C$ (or $\tl C'$) then $\tl{x_1}$ and
 $ \tl{x_2}$ are in the same coincidence class of $(\tl f,\tl g)$

\item $\tl C$ and $\tl C'$ are in the same coincidence class
of the pair $(\tl f,\tl g)$ if and only if $C$ is defective.

\end{enumerate}
\elem

\dem We need only to take $q:\tl N\to N$ as the identity map in the
Corollary~\ref{numcop}, Corollary~\ref{lifdef} and Lemma~\ref{lifpat}.
$\fimdem$

\bcor \label{corind} Under the hypotheses of lemma~\ref{or} we have:

\begin{enumerate}
\item If $C$ is not defective then  $\tl C$ and $\tl C'$ are two coincidence
  classes of the pair $(\tl f,\tl g)$ with $ind(\tl f,\tl g,\tl
  C)=-ind(\tl f,\tl g,\tl C')$  and  $|ind|(\tl f,\tl g,\tl C)=|ind(f,g,C)|$.

\item If $C$ is defective then $\tl C \cup \tl C'$ is a unique coincidence
  class of the pair $(\tl f,\tl g)$ with $ind(\tl f,\tl g,\tl
  C\cup \tl C')=0$.
\end{enumerate}
\ecor

\dem The first part is a consequence of lemma~\ref{or}, and the second can be
proved using the fact that $\tl M$ is the double orientable covering of $M$.
It is useful to remember that the pair $(\tl f,\tl g)$ is a pair of maps
between orientable manifolds. $\fimdem$

\bcor \label{lef} Under de hypotheses of lemma~\ref{or} we have:

\begin{enumerate}
\item  $L(\tl f,\tl g)=0$;

\item $N(\tl f, \tl g)$ is even;

\item $N(f,g)\ge \frac{ N(\tl f, \tl g)}{2}$;

\item If $N(\tl f, \tl g)=0$ then all coincidence classes with nonzero
  semi-index of the pair $(f,g)$ are defective.

\end{enumerate}
\ecor

\dem We need only to see that $p(Coin(\tl f,\tl g))=Coin(f,g)$, and that in the
pair $(\tl f,\tl g)$  we "lose" the defective classes.$\fimdem $

\section{Applications}

\bteo \label{apli} Let $f,g:M\to N$ be two maps between compact closed
manifolds of the same dimension such that $M$ is nonorientable and $N$ is
orientable. Suppose that 
$N$ is such that for all
orientable manifolds $M'$ of the same dimension of $N$ and all pairs of maps
$f',g':M'\to N$  we have $L(f',g')=0 \Rightarrow N(f',g')=0$. Then all
coincidence classes with nonzero semi-index of the pair  $(f,g)$ are
defective.  \eteo

\dem The hypotheses on $N$ are enough to show, using the notation of the proof
of theorem~\ref{or}, that $N(\tl f,\tl g)=0$. So by corollary~\ref{lef}, all
coincidence classes with nonzero semi-index of the pair  $(f,g)$ are
defective. $\fimdem$

We note that the hypotheses on the manifold $N$ in theorem~\ref{apli}, in
dimension greater then 2, are equivalent to the converse of Lefschetz
theorem. In dimension 2 these hypotheses are not equivalent but necessary for
the converse of Lefschetz theorem.

\bobs The following manifolds satisfy the hypotheses on the manifold $N$ in
the theorem~\ref{apli}:
\begin{enumerate}

\item Jiang spaces (\cite[Corollary 1]{dape1}).

\item Nilmanifolds (\cite[Theorem 5]{dape2}).

\item Homogeneous spaces of a compact connected Lie group $G$ by
  a finite subgroup $K$ (\cite[Theorem 4]{dape1}). 

\item Suitable manifolds (\cite[Theorem 3 e Theorem 4]{peter}) .

\end{enumerate}

\eobs

\section*{Acknowledgements}

This work was made during a postdoctoral year of the author at La-boratoire
\'Emile Picard - Universit\'e Paul Sabatier (Toulouse - France). We would like
to thank John Guaschi and Claude Hayat-Legrand for the invitation and
hospitality, and Peter N.-S. Wong for helpful conversations. This work was
supported by Capes-BEX0755/02-8 (international cooperation Capes-Cofecub
project number 364/01).

\end{document}